\documentclass[12pt]{amsart}
\usepackage{amsfonts}
\usepackage{amsthm}
\usepackage{graphics}
\usepackage{indentfirst}
\usepackage{cite}
\usepackage{latexsym}
\usepackage{amsmath}
\usepackage{amssymb}
\usepackage[dvips]{epsfig}
\usepackage{amscd}

\newtheorem{theorem}{Theorem}[section]
\newtheorem{remark}{Remark}[section]

\newtheorem{lemma}[theorem]{Lemma}

\newcommand{\bt}{\begin{theorem}}
\newcommand{\bl}{\begin{lemma}}
\newcommand{\el}{\end{lemma}}
\newcommand{\et}{\end{theorem}}
\newcommand{\br}{\begin{remark}}
\newcommand{\er}{\end{remark}}
\newcommand{\bc}{\begin{corollary}}
\newcommand{\ec}{\begin{corollary}}

\newcommand{\la}{\label}

\newcommand{\bn}{\begin{eqnarray}}
\newcommand{\en}{\end{eqnarray}}
\newcommand{\bnn}{\begin{eqnarray*}}
\newcommand{\enn}{\end{eqnarray*}}

\newcommand{\ba}{\begin{aligned}}
\newcommand{\ea}{\end{aligned}}
\newcommand{\be}{\begin{equation}}
\newcommand{\ee}{\end{equation}}

\def\norm[#1]#2{\|#2\|_{#1}}

\def\la{\label}

\begin{document}

\title{Global solvability of 3D inhomogeneous Navier-Stokes equations with density-dependent viscosity}

\author{Xiangdi Huang}
\address{NCMIS, Academy of Mathematics and Systems Science, CAS, Beijing 100190, P. R. China }
\email{xdhuang@amss.ac.cn}

\author{Yun Wang}
\address{Department of Mathematics, Soochow University, 1 Shizi Street, Suzhou 215006, P.R. China}
\email{ywang3@suda.edu.cn}
\date{}
\maketitle

\begin{abstract}In this paper, we consider the three-dimensional inhomogeneous Navier-Stokes
equations with density-dependent viscosity in presence of vacuum over bounded domains. Global-in-time unique strong solution is proved to exist when $\|\nabla u_0\|_{L^2}$ is suitably small with arbitrary large initial density. This generalizes all the previous results even for the constant viscosity.


AMS: 35Q35, 35B65, 76N10
\end{abstract}

\pagestyle{myheadings}
\thispagestyle{plain}
\markboth{XIANGDI HUANG AND YUN WANG}{3D DENSITY-DEPENDENT NAVIER-STOKES}

%%%%%%%%%%%%%%%Introduction%%%%%%%%%%%%%%%%%%%%%%%%%%%%%%%
\section{Introduction}

The Navier-Stokes equations are usually used to describe the motion of fluids. In particular,
for the study of multiphase fluids without surface tension, the following density-dependent Navier-Stokes equations acts as a model on some bounded domain $\Omega\subset R^N(N=2,3)$,
\be \la{ns} \left\{\ba
&\rho_t + {\rm div}~(\rho u) =0, \ \ \ \ \mbox{in}\ \Omega\times (0, T],\\
&(\rho u)_t + {\rm div}~(\rho u\otimes u) - {\rm div}~(2\mu(\rho) d) + \nabla P = 0,\ \ \ \ \mbox{in}\ \Omega\times (0, T],\\
&{\rm div}~ u =0,\ \ \ \ \ \ \mbox{in }\ \Omega \times [0, T],\\
&u=0,\ \ \ \ \ \ \ \mbox{on}\ \partial \Omega \times [0, T],\\
&\rho|_{t=0}= \rho_0, \ \ \ u|_{t=0} = u_0, \ \ \ \  \ \mbox{in }\ \Omega.
\ea \right. \ee
Here $\rho, u, $ and $P$ denote the density, velocity and pressure of the fluid, respectively. 
$$d = \frac12 \left[ \nabla u + (\nabla u)^T\right]$$
 is the deformation tensor.

$\mu = \mu(\rho)$ states the viscosity and is a function of $\rho$, which is assumed to satisfy
\be \la{mu}
\mu \in C^1[0, \infty),\ \  \mbox{and}\ \ \ \mu\geq \underline{\mu} >0\ \ \ \mbox{on}\ [0, \infty)\quad \mbox{for some positive constant $\underline{\mu}$}.
\ee
In this paper, we study the  initial boundary value problem to the system (\ref{ns})-(\ref{mu}) .

 The mathematical study for nonhomogeneous incompressible flow was initiated by the Russian school. They studied the case that
$\mu(\rho)$ is a constant and the initial density $\rho_0$ is bounded away from $0$. In the absence of vacuum, global existence of weak solutions as well as local strong solution was established by Kazhikov\cite{Kazhikov,AK}. The uniqueness of local strong solutions was first established by Ladyzhenskaya-Solonnikov\cite{LS}  for the initial boundary value problem, see also \cite{Salvi}. Furthermore, unique local strong solution is proved to be global in 2D\cite{Simon}. In recent years, Danchin initiated the studies for solutions in critical spaces. He\cite{Danchin-1, Danchin-2} derived the global well-posedness for small initial velocity in critical spaces, where density is close to a constant. For some subsequent works, refer to \cite{Abidi-Gui-Zhang, Paicu-Zhang} and references therein. We remark that in the very interesting papers \cite{Danchin-Mucha, Danchin-Mucha-2}, Danchin-Mucha studied the case for which density is piecewise constant, see also some generalizations in 2D \cite{Huang-Paicu-Zhang}.

When initial vacuum is taken into account and $\mu(\rho)$ is still a constant, Simon\cite{Simon} proved the global existence of weak solutions. Later, Choe-Kim\cite{Choe-Kim} proposed a compatibility condition as \eqref{compatibility} below to establish local existence of strong solution. Global strong solution allowing vacuum in 2D was recently derived by the authors\cite{HW}. Meanwhile,
some global solutions in 3D with small critical norms have been constructed, refer to the results in \cite{Abidi-Gui-Zhang-2, Craig-Huang-Wang} and references therein.

Finally, we come to the most general case: viscosity $\mu(\rho)$ depends on density $\rho$. Most results were concentrated on 2D case. Global weak solutions were derived by the revolutionary work \cite{DiPerna-Lions1, Lions} of DiPerna and Lions. Later,  Desjardins\cite{Desjardins} proved the global weak solution with more regularity for the two-dimensional case provided that the viscosity function $\mu(\rho)$ is a small pertubation of a  positive constant in $L^\infty$- norm.  Very recently, Abidi-Zhang\cite{Abidi-Zhang} generalized this 2D result to strong solutions.  Regarding the 3D case, Cho-Kim\cite{Cho-Kim} constructed a unique local strong solution by imposing some initial compatibility condition. Their result is stated as follows:

\bt\la{main-result-CK}
Assume that the initial data $(\rho_0, u_0)$ satisfies the regularity condition
\be \la{regularity} 0\leq \rho_0 \in W^{1, q}, \ \ 3< q < \infty , \ \ \ u_0\in H_{0, \sigma}^1 \cap H^2,
\ee
and the compatibility condition
\be\la{compatibility}
-{\rm div} \left( \mu(\rho_0) \left[  \nabla u_0 + (\nabla u_0)^T\right]\right) + \nabla P_0 = \rho_0^{\frac12} g~,
\ee
for some $(P_0, g)\in H^1 \times L^2$. Then there exists a small time $T $ and a unique strong solution $(\rho, u, P)$ to the initial boundary value problem \eqref{ns} such that
\be \nonumber \ba
& \rho \in C([0, T]; \ W^{1, q}),\ \ \ \ \ \nabla u, P \in C([0, T]; \ H^1) \cap L^2(0, T;\ W^{1,r}),\\
& \rho_t \in C ([0, T]; L^q), \ \ \ \ \ \ \sqrt{\rho} u_t \in L^\infty(0, T; L^2),\ \ \ u_t\in L^2(0, T; H_0^1),
\ea \ee
for any $r$ with $1\leq r < q$.
 Furthermore, if $T^*$ is the maximal existence time of the local strong solution $(\rho, u)$ , then either $T^*=\infty$ or
\be\la{criterion-CK}
\sup_{0\leq t < T^*} \left( \|\nabla \rho(t)\|_{L^q}   +  \|\nabla u(t)\|_{L^2}  \right)  = \infty.
\ee
\et

Motivated by the global existence result \cite{Craig-Huang-Wang} for the special case that $\mu$ is a constant, we aim to establish global well-posedness result for variable coefficient case. However, due to the strong coupling between viscosity coefficient and density, it's more complicated and involved with variable coefficient $\mu(\rho)$  and requires more delicate analysis.

Our main result proves the existence of global strong solution, provided $\|\nabla u_0\|_{L^2}$ is suitably small. To conclude, we arrive at
\bt\la{main-result}  Assume that the initial data $(\rho_0, u_0)$ satisfies \eqref{regularity}-\eqref{compatibility}, and $0 \leq \rho_0 \leq \overline{\rho}$.
Then there exists some small positive constant $\epsilon_0$, depending on $\Omega$, $q$, $\overline{\rho}$,    $\overline{\mu}= \displaystyle\sup_{[0, \overline{\rho}]} \mu(\rho)$, $\underline{\mu}$, $\|\nabla\mu(\rho_0)\|_{L^q}(=M)$, such that if
\be\label{mmain}
\|\nabla u_0\|_{L^2}\le \epsilon_0,
\ee
then the initial boundary value problem \eqref{ns} admits a unique global strong solution $(\rho, u)$, with
\be\la{global-solution-1} \ba
& \rho \in C([0, \infty) ; \ W^{1, q}),\ \ \ \ \ \nabla u, P \in C([0, \infty); \ H^1) \cap L^2_{loc}(0, \infty;\ W^{1,r}),\\
& \rho_t \in C ([0, \infty) ; L^q), \ \ \ \ \ \ \sqrt{\rho} u_t \in L^\infty_{loc}(0, \infty ; L^2),\ \ \ u_t\in L^2_{loc}(0, \infty; H_0^1),
\ea \ee
for any $r $, $1\leq r < q$.
\et

The main idea is to combine techniques developed by the authors in \cite{Huang-Li-Xin-CPAM, Craig-Huang-Wang} and time weighted energy estimates successfully applied to compressible Navier-Stokes equation  by Hoff\cite{Hoff}.

 Let's briefly sketch the proof. First we assume that $\|\nabla \mu(\rho)\|_{L^q}$ is less than $4M$ and $\|\nabla u\|_{L^2}^2$ is less than $4\|\nabla u_0\|_{L^2}^2$ on $[0, T]$,
then we prove that in fact $\|\nabla \mu(\rho)\|_{L^q}$ is less than $2M$ and $\|\nabla u\|_{L^2}^2$ is less than
$2 \|\nabla u_0\|_{L^2}^2$  on $[0, T]$,
under the assumption $\|\nabla u_0\|_{L^2}\leq \epsilon_0\leq \frac12$. On the other hand,
the control of $\|\nabla \mu(\rho)\|_{L^q}$ and $\|\nabla u\|_{L^2}$ lead to uniform estimates for other quantities, which guarantee the extension of local strong solutions. All the above procedures depends on a time independent bound of $\|\nabla u\|_{L^1L^\infty}$. Thanks to bounded domain, $u$ indeed has exponentially decay rather than insufficient polynomial decay for the whole space, that's the main reason why we can only treat system \eqref{ns} in bounded domain.

\begin{remark}
Theorem \ref{main-result} also hold for 2D case. Since the proof is quite similar, we omit it for simplicity. 
\end{remark}

The rest of the paper is organized as follows: Section 2 consists of some notations, definitions, and basic lemmas. Section 3 is devoted to the proof of Theorem \ref{main-result}.

%%%%%%%%%%%%%%%%%%%%%%%%%%%%Section 2%%%%%%%%%%%%%%%%%%%%%%%%%%%%%%%%%%%%%

\section{Preliminaries}

$\Omega$ is a bounded smooth domain in $\mathbb{R}^3$. Denote
$$ \int f \, dx = \int_\Omega f \,  dx.$$
For $1\leq r \leq \infty$ and $k\in \mathbb{N}$, the Sobolev spaces are defined in a standard way,
$$
L^r= L^r(\Omega), \ \ \ \ \ W^{k, r} = \left\{f\in L^r: \ \ \ \nabla^k f \in L^r \right\},\ \ $$
$$ H^k = W^{k, 2},\ \ \ \ C_{0, \sigma}^\infty = \{f\in C_0^\infty :\ \  {\rm div}~ f =0 \}.
$$
$$ H_0^1 = \overline{C_0^\infty},\ \ \ H_{0, \sigma}^1 = \overline{C_{0, \sigma}^\infty},\ \  \mbox{closure in the norm of } H^1.
$$

%%%%%%%%%%%%%%%%%%%%Lemma for Stokes%%%%%%%%%%%%%%%%%%%
High-order a priori estimates rely on the following regularity results for density-dependent Stokes equations.
\bl\la{stokes} Assume that $\rho \in W^{1, q}$, $3<q<\infty$, and $0 \leq \rho \leq  \bar{\rho}$.  Let $(u, P)\in H_{0, \sigma}^1 \times L^2$ be the unique weak solution to the boundary value problem
\be\la{stoke-system}
- {\rm div}~(2\mu (\rho) d) + \nabla P = F,\ \ \ \ {\rm div}~u =0\ \ \mbox{in}\ \Omega, \ \ \mbox{and}\ \ \int \frac{ P}{\mu(\rho)}\,  dx = 0,
\ee
where $d = \frac12\left[\nabla u + (\nabla u)^T\right]$ and
\be \nonumber
\mu \in C^1 [0, \infty),\ \ \ \underline{\mu} \leq \mu(\rho) \leq  \overline{\mu}\ \ \mbox{on} \ [0, \bar{\rho}].
\ee
 Then we have the following regularity results:

(1)\ If $F\in L^2$, then $(u, P)\in H^2 \times H^1$ and
\be\la{W22}
\|u\|_{H^2}  +  \|P/\mu(\rho) \|_{H^1}  \leq  C \left(  \frac{1}{\underline{\mu}} +     \frac{\overline{\mu}}{\underline{\mu}^{\frac{1}{\theta_2} +2 }}
\| \nabla \mu\|_{L^q}^{\frac{1}{\theta_2}}      \right) \|F\|_{L^2},
\ee
where $\theta_2$ satisfies
$$ \frac12 -\frac1q = \frac{\theta_2}{3} + \frac16,\ \ \ \ \mbox{i. e., }\ \ \ \theta_2 = \frac{q}{q-3}.$$

(2)\ If $F\in L^r$ for some $r\in (2, q)$, then $(u, P) \in W^{2, r} \times W^{1, r}$ and
\be \la{W2r}
\|u\|_{W^{2, r}} + \|P / \mu(\rho) \|_{W^{1,r}} \leq C \left(   \frac{1}{\underline{\mu}} +     \frac{\overline{\mu}}{\underline{\mu}^{\frac{1}{\theta_r} +2 }}
\| \nabla \mu\|_{L^q}^{\frac{1}{\theta_r}}                    \right) \|F\|_{L^r},
\ee
where
$$  \frac{1}{\theta_r} = \frac{\frac56 - \frac1r }{ \frac13 - \frac1q}.$$
Here the constant $C$ in \eqref{W22} and \eqref{W2r} depends on $\Omega$, $q$, $r$.
\el

The proof of Lemma \ref{stokes} has been given in \cite{Cho-Kim}, although the lemma is slightly different from the version in \cite{Cho-Kim}. We sketch it here for completeness.

\begin{proof}For the existence and uniqueness of the solution, please refer to Giaquinta-Modica\cite{GM}. We give the a priori estimates here. Assume that $F\in L^2$.
Multiply the first equation of \eqref{stoke-system} by $u$ and integrate over $\Omega$, then by Poincar\'e's inequality,
\be \nonumber
\int 2 \mu(\rho) |d|^2 \, dx = \int F\cdot u \, dx \leq \|F\|_{L^2} \cdot \|u\|_{L^2}\leq C \|F\|_{L^2} \cdot \|\nabla u\|_{L^2}.
\ee
Note that $$\displaystyle 2\int |d|^2 \, dx = \int |\nabla u|^2 \, dx , $$ hence
\be \label{2.2}
\|\nabla u\|_{L^2} \leq C \underline{\mu}^{-1} \|F\|_{L^2}.
\ee

Since $\displaystyle \int \frac{P}{\mu(\rho)} \, dx =0$, according to Bovosgii's theory, there exists a function $v\in H_0^1$, such that
\be \nonumber
 {\rm div}~v =\frac{P}{\mu(\rho)},
\ee
and
\be \nonumber
\|\nabla v\|_{L^2} \leq C \left\|  \frac{P}{\mu(\rho)}  \right\|_{L^2}.
\ee
Multiply the first equation of \eqref{stoke-system} by $-v$, and integrate over $\Omega$, then
\be \nonumber \ba
\int \frac{P^2}{\mu(\rho) } \, dx &= -\int F\cdot v \, dx + 2 \int \mu(\rho) d: \nabla v \, dx\\
& \leq \|F\|_{L^2} \cdot \|v\|_{L^2} + C\overline{\mu} \cdot  \|\nabla u\|_{L^2} \cdot \|\nabla v\|_{L^2} \\
& \leq C \|F\|_{L^2} \cdot \|\nabla v\|_{L^2} + C \frac{\overline{\mu}}{\underline{\mu}} \cdot  \|F\|_{L^2} \cdot \|\nabla v\|_{L^2} \\
& \leq C \frac{\overline{\mu}}{\underline{\mu}}  \cdot \|F\|_{L^2} \cdot \left\| \frac{P}{\mu(\rho)}  \right\|_{L^2}.
\ea \ee
On the other hand side,
\be \nonumber
\int \frac{P^2}{\mu(\rho)} \, dx \geq \underline{\mu} \int \frac{P^2}{\mu(\rho)^2} \, dx.
\ee
Hence,
\be \la{2.3}
\left\| \frac{P}{\mu(\rho)}   \right\|_{L^2} \leq C \frac{\overline{\mu}}{\underline{\mu}^2}  \cdot \|F\|_{L^2} .
\ee

The first equation of \eqref{stoke-system} can be re-written as
\be \nonumber
-\Delta u + \nabla \left( \frac{P}{\mu(\rho)}  \right) = \frac{F}{\mu(\rho)}
+ \frac{2d\cdot \nabla \mu(\rho)}{\mu(\rho)} - \frac{P \nabla \mu(\rho) }{ \mu(\rho)^2}.
\ee
By virture of the classical theory for Stokes equations and Gagliardo-Nirenberg inequality, we have
\be \nonumber  \ba
& \|u\|_{H^2} + \left\| \nabla \left(  \frac{P}{\mu(\rho) }  \right) \right\|_{L^2} \\
& \leq C \left(    \left\| \frac{F}{\mu(\rho)}   \right\|_{L^2}     + \left\|  \frac{d \cdot \nabla \mu(\rho)}{\mu(\rho)} \right\|_{L^2}
+ \left\|  \frac{P \nabla \mu(\rho) }{ \mu(\rho)^2} \right\|_{L^2} \right)\\
& \leq C \left( \underline{\mu}^{-1} \|F\|_{L^2} + \underline{\mu}^{-1} \|\nabla \mu\|_{L^q}    \cdot
\|\nabla u \|_{L^{\frac{2q}{q-2}}}   + \underline{\mu}^{-1} \|\nabla \mu\|_{L^q} \cdot \left\|   \frac{P}{\mu}\right\|_{L^{\frac{2q}{q-2}}}   \right) \\
& \leq C \left[  \underline{\mu}^{-1} \|F\|_{L^2} + \underline{\mu}^{-1} \|\nabla \mu\|_{L^q} \cdot \|\nabla u\|_{L^2}^{\theta_2} \cdot \|\nabla u\|_{H^1}^{1- \theta_2} \right. \\
&\ \ \ \ \ \ \ \ \ \ \ \ \ + \left. \underline{\mu}^{-1} \|\nabla \mu\|_{L^q} \cdot \left\| \frac{P}{\mu}  \right\|_{L^2}^{\theta_2} \cdot  \left\|  \nabla \left( \frac{P}{\mu}   \right)   \right\|_{L^2}^{1-\theta_2} \right].
\ea \ee
By Young's inequality,
\be \la{2.6} \ba
& \|u\|_{H^2} + \left\| \nabla \left(  \frac{P}{\mu(\rho) }  \right) \right\|_{L^2} \\
& \leq C \underline{\mu}^{-1} \|F\|_{L^2} + C \underline{\mu}^{-\frac{1}{\theta_2}} \|\nabla \mu\|_{L^q}^{\frac{1}{\theta_2}} \cdot \left(  \|\nabla u\|_{L^2} + \left\|   \frac{P}{\mu} \right\|_{L^2}\right)\\
& \leq C \underline{\mu}^{-1} \|F\|_{L^2} + C \underline{\mu}^{-\frac{1}{\theta_2} -1 } \cdot \left( 1+ \frac{\overline{\mu}}{\underline{\mu}}   \right) \cdot \|\nabla \mu\|_{L^q}^{\frac{1}{\theta_2}} \cdot \|F\|_{L^2} \\
& \leq C \left(  \frac{1}{\underline{\mu}} +     \frac{\overline{\mu}}{\underline{\mu}^{\frac{1}{\theta_2} +2 }}
\| \nabla \mu\|_{L^q}^{\frac{1}{\theta_2}}      \right) \|F\|_{L^2},
\ea
\ee
where $\theta_2$ satisfies
$$ \frac12 -\frac1q = \frac{\theta_2}{3} + \frac16,\ \ \ \ \mbox{or}\ \ \ \theta_2 = \frac{q}{q-3}.$$

Similarly,
\be \la{2.7}
\| u\|_{W^{2, r}} + \left\|  \nabla \left( P/\mu(\rho)   \right)  \right\|_{L^r}
\leq C \left(   \frac{1}{\underline{\mu}} +     \frac{\overline{\mu}}{\underline{\mu}^{\frac{1}{\theta_r} +2 }}
\| \nabla \mu\|_{L^q}^{\frac{1}{\theta_r}}                    \right) \|F\|_{L^r},
\ee
where
$$  \frac{1}{\theta_r} = \frac{\frac56 - \frac1r }{ \frac13 - \frac1q}.$$
\end{proof}

%%%%%%%%%%%%%%%%%%%%%%%%%%%%%%%Section3%%%%%%%%%%%%%%%%%%%%%%%%%%%%%%%%%%%%%

\section{Proof of Theorem \ref{main-result}}
The proof of Theorem \ref{main-result} is composed of two parts. The first part contains a priori time-weighted estimates of different levels. Upon these estimates, the second part uses a contradiction induction process to extend the local strong solution. The two parts are presented in Subsections 3.1 and 3.2, respectively.

\subsection{A Priori Estimates}
In this subsection, we establish some a priori time-weighted estimates. The initial velocity belongs to $H^1$, but some uniform estimates of higher order and independent of time are required.
To achieve that, we take some power of time as a weight. The idea is based on the parabolic property of the system. In this subsection, the constant $C$ will denote some positive constant which maybe dependent on $\Omega$, $q$, but is independent of $\rho_0$ or $u_0$.

%%%%%%%%%%%%%%%%%%%Lemma for density%%%%%%%%%%%%%%%%
First, as the density satisfies the transport equation $(\ref{ns})_1$ and making use of $(\ref{ns})_3$, one has the following lemma.
\bl\la{density-3}Suppose $(\rho, u, P)$ is the unique local strong solution to \eqref{ns} on $[0, T]$, with the initial data $(\rho_0, u_0)$, it holds that
\be \nonumber
0 \leq \rho(x, t) \leq \bar{\rho}, \ \ \ \ \ \mbox{for every}\ (x, t)\in \Omega \times [0, T].
\ee
\el

%%%%%%%%%%%%%%%%%%Lemma for Basic Energy Inequality%%%%%%%%%%%%%
Next, the basic energy inequality of the system (\ref{ns}) reads
\bt \la{energy-inequality-3}Suppose $(\rho, u, P)$ is the unique local strong solution to \eqref{ns} on $[0, T]$, with the initial data $(\rho_0, u_0)$, it holds that
\be\la{energyinequality3-1}
\int \rho |u(t)|^2 \, dx + \int_0^t \int \mu(\rho) |d|^2 \, dx ds \leq C \cdot \overline{\rho}\cdot\| u_0 \|_{L^2}^2 \, , \ \ \ \mbox{for every } \ t\in [0, T],
\ee
or in other words,
\be\la{energyinequality3}
\int \rho |u(t)|^2 \, dx + \underline{\mu} \int_0^t \int |\nabla u|^2 \, dx ds \leq C \cdot \overline{\rho}\cdot\| u_0 \|_{L^2}^2 \, , \ \ \ \mbox{for every } \ t\in [0, T].
\ee
\et

\begin{proof}The proof  is standard. Multiplying the momentum equation by $u$ and integrating over $\Omega$ yield that
\be \nonumber
\frac12 \frac{d}{dt} \int \rho |u|^2 \,  dx + 2\int \mu(\rho) |d|^2 \, dx = 0~.
\ee
Then \eqref{energyinequality3} is true owing to the fact $2 \int |d|^2 dx = \int |\nabla u|^2 dx$ and $\mu(\rho) \geq \underline{\mu}$.
\end{proof}

%%%%%%%%%%%%%%%%%%%%%%%%%%Prop1 in Section 4%%%%%%%%%%%%%%%%%%%%%%%%%%%%%%%%%%%
Denote
\be \nonumber
M = \|\nabla \mu(\rho_0)\|_{L^q},
\ee
\be \nonumber
 M_2 = \frac{1}{\underline{\mu}} + \frac{\overline{\mu}}{\underline{\mu}} \cdot \frac{1}{\underline{\mu}^{1 /\theta_2 + 1 }} \cdot (4M)^{\frac{1}{\theta_2}},
\ee
and
\be \nonumber M_r = \frac{1}{\underline{\mu}} + \frac{\overline{\mu}}{\underline{\mu}} \cdot \frac{1}{\underline{\mu}^{1 /\theta_r + 1 }} \cdot (4M)^{\frac{1}{\theta_r}}.
\ee

\bt \label{prop1-section4}
 Suppose $(\rho, u, P)$ is the unique local strong solution to \eqref{ns} on $[0, T]$, with the initial data $(\rho_0, u_0)$, and satisfies
\be\la{condition1}
\sup_{t\in [0, T]} \|\nabla \mu(\rho(t))\|_{L^q} \leq 4M,
\ee
and \be \la{condition2}
\sup_{t\in [0, T]} \|\nabla u (t) \|_{L^2}^2  \leq 4\|\nabla u_0\|_{L^2}^2 \leq 1.
\ee
There exists a positive number $C_1$, depending on $\Omega$, $q$ such that if
\be \la{condition3}
 C_1 \underline{\mu}^{-2} (M_2^2 + M^4 M_2^6)\overline{\rho}^4  \cdot \|\nabla u_0\|_{L^2}^2  \leq \ln 2,
\ee
then
\be \la{3.2}
 \frac{1}{\underline{\mu}} \int_0^T \|\sqrt{\rho} u_t\|_{L^2}^2 \, dt+ \sup_{t\in [0, T]} \|\nabla u(t)\|_{L^2}^2 \leq 2\|\nabla u_0\|_{L^2}^2.
\ee
\et

Before the proof of Theorem \ref{prop1-section4}, let us introduce an auxiliary lemma, which is a result of the $W^{2, 2}$-estimates in Lemma \ref{stokes}.

%%%%%%%%%%%%%%%%%%%%%%%%Lemma for W^{2,2}%%%%%%%%%%%%%%%%%%%%%%
\bl \la{Alemma}  Suppose $(\rho, u, P)$ is the unique local strong solution to \eqref{ns} on $[0, T]$ and satisfies
\be\nonumber
\sup_{t\in [0, T]} \|\nabla \mu(\rho(t))\|_{L^q} \leq 4M.
\ee
Then it holds that
\be \nonumber
\|\nabla u\|_{H^1} \leq C M_2 \|\rho u_t\|_{L^2} + C M_2^2\cdot \overline{\rho}^2 \cdot \|\nabla u\|_{L^2}^3.
 \ee \el

\begin{proof} The momentum equation can be rewritten as follows,
\be \la{3.4}
-2 {\rm div}~(\mu(\rho) d ) + \nabla P = - \rho u_t - (\rho u \cdot \nabla) u.
\ee
It follows from Lemma \ref{stokes} and Gagliardo-Nirenberg inequality that
\be \nonumber \ba
\|\nabla u\|_{H^1} & \leq C M_2 \left(   \|\rho u_t\|_{L^2} + \|\rho u\cdot \nabla u\|_{L^2}            \right) \\
& \leq C M_2 \|\rho u_t\|_{L^2} + C M_2 \cdot \overline{\rho} \cdot \|u\|_{L^6} \cdot \|\nabla u\|_{L^3} \\
& \leq C M_2 \|\rho u_t\|_{L^2} + C M_2 \cdot \overline{\rho} \cdot \|\nabla u\|_{L^2}^{\frac32} \cdot \|\nabla u\|_{H^1}^{\frac12} .
\ea \ee
By Young's inequality,
\be \nonumber
\|\nabla u\|_{H^1} \leq C M_2 \|\rho u_t\|_{L^2} + C M_2^2 \cdot \overline{\rho}^2 \cdot \|\nabla u\|_{L^2}^3.
\ee
\end{proof}

\begin{proof}[Proof of Theorem \ref{prop1-section4}]
Multiply the momentum equation by $u_t$ and integrate over $\Omega$, then
\be\nonumber \ba
& \int \rho |u_t|^2 \, dx + \frac{d}{dt} \int \mu(\rho) |d|^2 \, dx \\
& \leq \left| \int \rho u \cdot \nabla u \cdot u_t\, dx        \right| + C \int |\nabla \mu(\rho) |\cdot |u|\cdot |\nabla u|^2 \, dx .
\ea
\ee
Here we used the renormalized mass equation for $\mu(\rho)$,
\be \nonumber
\partial_t [\mu(\rho)] + u \cdot \nabla \mu(\rho) = 0,
\ee
which is derived due to the fact ${\rm div}~u =0$.

Applying Gagliardo-Nirenberg inequality and Lemma \ref{Alemma},
\be\nonumber \ba
& \left| \int \rho u \cdot \nabla u \cdot u_t\, dx        \right| \\
& \leq  \frac18 \|\sqrt{\rho}  u_t\|_{L^2}^2  + C \|\sqrt{\rho} u\|_{L^6}^2 \cdot \|\nabla u\|_{L^3}^2\\
& \leq  \frac18 \|\sqrt{\rho}  u_t\|_{L^2}^2 + C \overline{\rho} \cdot \|\nabla u\|_{L^2}^3 \cdot \|\nabla u\|_{H^1}\\
& \leq  \frac18 \|\sqrt{\rho}  u_t\|_{L^2}^2 + C \overline{\rho} \cdot \|\nabla u\|_{L^2}^3 \cdot \left[  C M_2 \|\rho u_t\|_{L^2} + C M_2^2 \overline{\rho}^2 \|\nabla u\|_{L^2}^3    \right],
\ea
\ee
and similarly,
\be \nonumber  \ba
&  C \int |\nabla \mu(\rho) |\cdot |u|\cdot |\nabla u|^2 \, dx \\
& \leq  C \|\nabla \mu(\rho) \|_{L^3} \cdot \|u\|_{L^6} \cdot \|\nabla u\|_{L^4}^2 \\
& \leq  C \|\nabla \mu(\rho) \|_{L^3} \cdot \|\nabla u\|_{L^2}^{\frac32} \cdot \|\nabla u\|_{H^1}^{\frac32}\\
& \leq  C M \|\nabla u\|_{L^2}^{\frac32} \cdot \left[   C M_2 \|\rho u_t\|_{L^2} + C M_2^2 \overline{\rho}^2 \|\nabla u\|_{L^2}^3      \right]^{\frac32}.
\ea \ee
Hence, by  Young's inequality,
\be \nonumber  \ba
& \int \rho |u_t|^2 \, dx + \frac{d}{dt} \int \mu(\rho) |d|^2 \, dx \\
& \leq \frac18 \|\sqrt{\rho}  u_t\|_{L^2}^2 + \frac18 \|\sqrt{\rho} u_t\|_{L^2}^2 + C  M_2^2 \overline{\rho}^3 \|\nabla u\|_{L^2}^6 +  \frac18 \|\sqrt{\rho}  u_t\|_{L^2}^2 \\
&\ \ \ \ \  + C \left(  M M_2^{\frac32} \cdot \overline{\rho}^{\frac34} \|\nabla u\|_{L^2}^{\frac32}  \right)^4
+ C M M_2^3 \cdot \overline{\rho}^3 \cdot \|\nabla u\|_{L^2}^6\\
& \leq \frac38 \|\sqrt{\rho}  u_t\|_{L^2}^2 +  C \left( M_2^2 +  M^4 M_2^6 + M M_2^3  \right) \cdot \overline{\rho}^3 \cdot \|\nabla u\|_{L^2}^6.
\ea \ee
So we have
\be \la{3.6}
 \int \rho |u_t|^2 \, dx + \frac{d}{dt} \int \mu(\rho) |d|^2 \, dx \\
\leq C \left( M_2^2 +  M^4 M_2^6  \right) \cdot \overline{\rho}^3 \cdot \|\nabla u\|_{L^2}^6.
\ee
Integrate with respect to time on $[0, t]$,
\be \nonumber
\frac{1}{\underline{\mu}} \int_0^t \int \rho |u_t|^2 \, dx ds + \sup_{s\in [0, t]} \int |\nabla u|^2 \, dx
\leq C \underline{\mu}^{-1} \left( M_2^2 +  M^4 M_2^6  \right) \cdot \overline{\rho}^3
\cdot \int_0^t \|\nabla u\|_{L^2}^6 \, ds.
\ee
Applying Gronwall's inequality,
\be \nonumber \ba
& \frac{1}{\underline{\mu}} \int_0^T \int \rho |u_t|^2 \, dx dt + \sup_{t\in [0, T] } \int |\nabla u|^2 \, dx \\
& \leq \|\nabla u_0\|_{L^2}^2 \cdot \exp\left\{ C \underline{\mu}^{-1} (M_2^2 + M^4 M_2^6 ) \cdot \overline{\rho}^3 \cdot \int_0^T \|\nabla u\|_{L^2}^4 \, dt \right\} .
\ea \ee
According to Theorem \ref{energy-inequality-3} and the assumption \eqref{condition2},
\be \la{3.7} \ba
\int_0^T \|\nabla u\|_{L^2}^4 \, dt & \leq \sup_{t\in [0, T]} \|\nabla u\|_{L^2}^2 \cdot \int_0^T \|\nabla u\|_{L^2}^2 \, ds \\
& \leq C \underline{\mu}^{-1}\cdot \overline{\rho} \cdot \|u_0\|_{L^2}^2 \\
& \leq C \underline{\mu}^{-1}\cdot \overline{\rho}\cdot \|\nabla u_0\|_{L^2}^2.
 \ea \ee
Hence, we arrive at
\be \la{3.8} \ba
& \frac{1}{\underline{\mu}} \int_0^T \int \rho |u_t|^2 \, dx dt + \sup_{t\in [0, T] } \int |\nabla u|^2 \, dx \\
& \leq  \|\nabla u_0\|_{L^2}^2 \exp\{ C_1 \underline{\mu}^{-2} (M_2^2 + M^4 M_2^6) \cdot \overline{\rho}^4
\cdot \|\nabla u_0\|_{L^2}^2 \}.
\ea \ee
Now it is clear that  \eqref{3.2} holds, provided \eqref{condition3} holds.

\end{proof}

As a byproduct of the estimates in the proof, we have the following result.
\bt \la{corollary}
Suppose $(\rho, u, P)$ is the unique local strong solution to \eqref{ns} on $[0, T]$, with the initial data $(\rho_0, u_0)$,  and satisfies
the assumptions \eqref{condition1}-\eqref{condition3} as in Theorem \ref{prop1-section4}.
Then
\be \la{3.10}
\frac{1}{\underline{\mu}} \int_0^T t \|\sqrt{\rho} u_t\|_{L^2}^2 \, dt  + \sup_{t\in [0, T]} t\|\nabla u\|_{L^2}^2
  \leq \frac{C \cdot  \overline{\rho}}{\underline{\mu} } \|\nabla u_0\|_{L^2}^2.
\ee
\et

\begin{proof}Multiplying \eqref{3.6} by $t$, as shown in the last proof, one has
\be \la{3.11} \ba
& \frac{1}{\underline{\mu}} \int_0^T t \|\sqrt{\rho} u_t\|_{L^2}^2 \, dt
+ \sup_{t\in [0, T]} t\|\nabla u\|_{L^2}^2 \\
& \leq \frac{1}{\underline{\mu}} \int_0^T  \int \mu(\rho) |d|^2 \, dx  dt \cdot  \exp\{ C_1 \underline{\mu}^{-2} (M_2^2 + M^4 M_2^6) \cdot \overline{\rho}^4 \cdot \|\nabla u_0\|_{L^2}^2 \}.
\ea \ee
According to Theorem \ref{energy-inequality-3},
\be \la{3.12}
\int_0^T  \int \mu(\rho) |d|^2 \, dx dt \leq C \cdot  \overline{\rho} \cdot \|u_0\|_{L^2}^2
\leq  C \cdot \overline{\rho} \cdot \|\nabla u_0\|_{L^2}^2.
 \ee
Hence,
\be\nonumber \ba
& \frac{1}{\underline{\mu}} \int_0^T t \|\sqrt{\rho} u_t\|_{L^2}^2 \, dt  + \sup_{t\in [0, T]} t\|\nabla u\|_{L^2}^2 \\
& \leq \frac{C\overline{\rho}}{\underline{\mu}}\|\nabla u_0\|_{L^2}^2 \cdot
\exp\{ C_1 \underline{\mu}^{-2} (M_2^2 + M^2 M_2^6 ) \cdot \overline{\rho}^4 \cdot \|\nabla u_0\|_{L^2}^2    \} \\
& \leq \frac{C\cdot \overline{\rho}}{\underline{\mu}} \|\nabla u_0\|_{L^2}^2.
\ea \ee

\end{proof}

%%%%%%%%%%%%%%%%%%%%%%%%%%Prop2 in Section 4%%%%%%%%%%%%%%%%%%%%%%%%%%%%%%%%%%

\bt
\label{prop2-section4}
 Suppose $(\rho, u, P)$ is the unique local strong solution to \eqref{ns} on $[0, T]$, with the initial data $(\rho_0, u_0)$,  and satisfies
the assumptions \eqref{condition1}-\eqref{condition3}.
Then
\be \la{4-5}  \ba
&\sup_{t\in [0, T] } t\int \rho |u_t|^2\, dx +\underline{\mu} \int_0^T t \|\nabla u_t\|_{L^2}^2 \\
&\leq C\|\nabla u_0\|_{L^2}^2\Theta_1 \cdot \exp\{C\Theta_2\}
\ea\ee
and
\be \la{4-5-add}  \ba
&\sup_{t^2\in [0, T] } t^2\int \rho |u_t|^2\, dx +\underline{\mu} \int_0^T t^2 \|\nabla u_t\|_{L^2}^2 \\
&\leq C\frac{\overline{\rho}}{\underline{\mu}}\|\nabla u_0\|_{L^2}^2\Theta_1 \cdot \exp\{C\Theta_2\}.
\ea\ee
where
\be
\Theta_1 = \frac{ M_2^4 \overline{\rho}^8}{\underline{\mu}^3} + \frac{ M^2 M_2^8 \overline{\rho}^{10}}{\underline{\mu}^3} +  \underline{\mu},\quad \Theta_2 =\frac{\overline{\rho}^4 }{\underline{\mu}^4 } + \frac{ M_2^2 \overline{\rho}^4 }{\underline{\mu}^2 } +  M^2 M_2^4 \overline{\rho}^2.
\ee
\et

\begin{proof}Take t-derivative of the momentum equation,
\be \la{4-5-0}
\rho u_{tt} + (\rho u) \cdot \nabla u_t - {\rm div}~(2\mu(\rho) d_t) + \nabla P_t = -\rho_t u_t
-(\rho u)_t \cdot \nabla u + {\rm div}~(2\mu(\rho)_t d).
\ee
Multiplying \eqref{4-5-0} by $t u_t$ and integrating over $\Omega$, we get after integration by parts that
\be \la{4-5-1}\ba
& \frac{t}{2} \frac{d}{dt} \int \rho |u_t|^2 \, dx +  2t \int \mu(\rho)  |d_t|^2 \, dx\\
&  = -t \int \rho_t |u_t|^2 \, dx - t\int (\rho u)_t \cdot \nabla u \cdot u_t \, dx
-t \int 2\mu(\rho)_t  d \cdot \nabla u_t \, dx~.
\ea \ee
Let us estimate the terms on the righthand of \eqref{4-5-1}. First, utilizing the mass equation and Poincar\'e's inequality, one has
\be \la{4-5-2}\ba
& -t\int \rho_t |u_t|^2\,  dx \\ &= - 2t\int \rho u \cdot \nabla u_t \cdot u_t \, dx\\
& \leq C\overline{\rho}^{\frac12} \cdot t \cdot \|\sqrt{\rho} u_t\|_{L^3}\cdot \|\nabla u_t\|_{L^2} \cdot \|u\|_{L^6}\\
& \leq C\overline{\rho}^{\frac12} \cdot t \cdot \|\sqrt{\rho} u_t\|_{L^2}^{\frac12} \cdot \|\sqrt{\rho} u_t\|_{L^6}^{\frac12} \cdot \|\nabla u_t\|_{L^2} \cdot \|\nabla u\|_{L^2}\\
& \leq C \overline{\rho}^{\frac34} \cdot t \cdot \|\sqrt{\rho} u_t\|_{L^2}^{\frac12} \cdot \|\nabla u_t\|_{L^2}^{\frac32}
\cdot \|\nabla u\|_{L^2} \\
& \leq \frac18 \underline{\mu}\cdot t \|\nabla u_t\|_{L^2}^2 + C \underline{\mu}^{-3} \cdot \overline{\rho}^3\cdot t\|\sqrt{\rho}u_t\|_{L^2}^2 \cdot \|\nabla u\|_{L^2}^4.
\ea
\ee

Second, utilizing the renormalized mass equation for $\mu(\rho)$,
\be \la{4-5-3}\ba
& -2t\int  \mu(\rho)_t \cdot d\cdot \nabla u_t\, dx \\
& \leq Ct  \int |u| \cdot |\nabla \mu(\rho) | \cdot |d| \cdot |\nabla u_t| \, dx \\
& \leq Ct \cdot  \|\nabla \mu(\rho) \|_{L^3} \cdot \|\nabla u_t\|_{L^2} \cdot \|d\|_{L^6} \cdot \|u\|_{L^\infty}\\
& \leq C M t \cdot \|\nabla u_t\|_{L^2} \cdot \|\nabla u\|_{H^1}^2 .
\ea \ee
It follows from Lemma \ref{Alemma} that
\be \la{4-5-4}  \ba
& - 2t \int  \mu(\rho)_t \cdot d\cdot \nabla u_t\, dx \\
& \leq \frac18 \underline{\mu} \cdot t\|\nabla u_t\|_{L^2}^2
+ \frac{CM^2}{\underline{\mu}}  \cdot t  \left(  M_2^4 \overline{\rho}^2 \|\sqrt{\rho} u_t\|_{L^2}^4
+ M_2^8 \overline{\rho}^8 \|\nabla u\|_{L^2}^{12}  \right)\\
& \leq \frac18 \underline{\mu} \cdot t\|\nabla u_t\|_{L^2}^2 + \frac{C M^2 M_2^4 \overline{\rho}^2 }{\underline{\mu}} \cdot t \|\sqrt{\rho} u_t\|_{L^2}^4 + \frac{C M^2 M_2^8 \overline{\rho}^8}{\underline{\mu}} \cdot t\|\nabla u\|_{L^2}^{12} .
\ea \ee

Finally, taking into account the mass equation again, we arrive at
\be \la{4-5-5} \ba
& - t\int (\rho u)_t \cdot \nabla u \cdot u_t \, dx \\
& =- t\int \rho u\cdot \nabla ( u\cdot \nabla u\cdot u_t ) \, dx - t\int \rho u_t \cdot \nabla u\cdot u_t \, dx \\
& \leq t\int \rho |u|\cdot |\nabla u|^2\cdot |u_t|  \, dx + Ct \int \rho |u|^2 \cdot |\nabla^2 u| \cdot |u_t| \, dx \\
&\ \ \ \ +t \int  \rho |u|^2 \cdot |\nabla u| \cdot |\nabla u_t| \, dx + t\int \rho |u_t|^2 \cdot |\nabla u|\, dx\\
& = \sum_{i=1}^4 J_i.
\ea \ee
Hence, it follows from Sobolev embedding inequality, Gagliardo-Nirenberg inequality, and Lemma \ref{Alemma} that
\be \la{4-5-6}
\ba
J_1 & \leq C \overline{\rho} \cdot t\|u_t\|_{L^6} \cdot \|u\|_{L^6} \cdot \|\nabla u\|_{L^3}^2\\
& \leq C \overline{\rho} \cdot t \|\nabla u_t\|_{L^2} \cdot \|\nabla u\|_{L^2}^2 \cdot \|\nabla u\|_{H^1}\\
& \leq C \overline{\rho} \cdot t \|\nabla u_t\|_{L^2} \cdot \|\nabla u\|_{L^2}^2 \cdot \left( M_2 \overline{\rho}^{\frac12}  \|\sqrt{\rho} u_t\|_{L^2}  + M_2^2 \overline{\rho}^2 \|\nabla u\|_{L^2}^3 \right) \\
& \leq \frac18 \underline{\mu} \cdot t\|\nabla u_t\|_{L^2}^2 + \frac{CM_2^2 \overline{\rho}^3  }{\underline{\mu}}
\cdot t\|\sqrt{\rho} u_t\|_{L^2}^2 \cdot \|\nabla u\|_{L^2}^4 +
\frac{CM_2^4 \overline{\rho}^6 }{\underline{\mu}}\cdot t \|\nabla u\|_{L^2}^{10}.
\ea
\ee
Similarly, it holds that
\be \la{4-5-7} \ba
J_2 & \leq C\overline{\rho} \cdot  t \|u_t\|_{L^6} \cdot \|\nabla^2 u\|_{L^2} \cdot \|u\|_{L^6}^2\\
& \leq C\overline{\rho} \cdot t  \| \nabla u_t\|_{L^2} \cdot \|\nabla u\|_{H^1} \cdot \|\nabla u\|_{L^2}^2\\
& \leq  \frac18 \underline{\mu} \cdot t\|\nabla u_t\|_{L^2}^2 + \frac{CM_2^2 \overline{\rho}^3  }{\underline{\mu}}
\cdot t\|\sqrt{\rho} u_t\|_{L^2}^2 \cdot \|\nabla u\|_{L^2}^4 +
\frac{CM_2^4 \overline{\rho}^6 }{\underline{\mu}}\cdot t \|\nabla u\|_{L^2}^{10},
\ea \ee
and
\be \la{4-5-8}
\ba
J_3 & \leq C \overline{\rho} \cdot t \|\nabla u_t\|_{L^2} \cdot \|\nabla u\|_{L^6} \cdot \|u\|_{L^6}^2 \\
& \leq  \frac18 \underline{\mu} \cdot t\|\nabla u_t\|_{L^2}^2 + \frac{CM_2^2 \overline{\rho}^3  }{\underline{\mu}}
\cdot t\|\sqrt{\rho} u_t\|_{L^2}^2 \cdot \|\nabla u\|_{L^2}^4 +
\frac{CM_2^4 \overline{\rho}^6 }{\underline{\mu}}\cdot t \|\nabla u\|_{L^2}^{10}.
\ea
\ee

Owing to Lemma \ref{Alemma} and Sobolev embedding inequality,
\be \la{4-5-9} \ba
J_4&  \leq C t \|\sqrt{\rho} u_t\|_{L^4}^2 \cdot \|\nabla u\|_{L^2} \\
& \leq C  t \|\sqrt{\rho} u_t\|_{L^2}^{\frac12}  \cdot \overline{\rho}^{\frac34} \|\nabla u_t\|_{L^2}^{\frac32} \cdot
\|\nabla u\|_{L^2}\\
& \leq \frac18 \underline{\mu}\cdot t \|\nabla u_t\|_{L^2}^2 + \frac{C \overline{\rho}^3}{\underline{\mu}} \cdot t\|\sqrt{\rho} u_t\|_{L^2}^2 \cdot \|\nabla u\|_{L^2}^4.
 \ea \ee

Combine all the above estimates \eqref{4-5-2}-\eqref{4-5-9},
\be \la{4-5-10} \ba
 &   \frac{d}{dt} t \int \rho |u_t|^2 \, dx +  \underline{\mu} \cdot t \|\nabla u_t\|_{L^2}^2 \\
& \leq \frac{C \overline{\rho}^3 }{\underline{\mu}^3} \cdot t\|\sqrt{\rho} u_t\|_{L^2}^2 \cdot \|\nabla u\|_{L^2}^4
                  + \frac{C M_2^2 \overline{\rho}^3}{\underline{\mu}}  \cdot t\|\sqrt{\rho} u_t\|_{L^2}^2 \cdot \|\nabla u\|_{L^2}^4  \\ &\ \ \ \ + \frac{C M_2^4 \overline{\rho}^6}{\underline{\mu}} \cdot t\|\nabla u\|_{L^2}^{10}
+ \frac{C M^2 M_2^4 \overline{\rho}^2 }{\underline{\mu}} \cdot t\|\sqrt{\rho} u_t\|_{L^2}^4 \\
&\ \ \ \ + \frac{C M^2 M_2^8 \overline{\rho}^8}{\underline{\mu}}\cdot t\|\nabla u\|_{L^2}^{12}
+ \int \rho |u_t|^2 \, dx.
\ea \ee
Applying Gronwall's inequality,
\be\nonumber \ba
& \sup_{t\in [0, T] } t\int \rho |u_t|^2\, dx +\underline{\mu} \int_0^T t \|\nabla u_t\|_{L^2}^2\, dt  \\
&  \leq \left[  \int_0^T  \left(  \frac{C M_2^4 \overline{\rho}^6}{\underline{\mu}}\cdot t\|\nabla u\|_{L^2}^{10} +
\frac{CM^2 M_2^8 \overline{\rho}^8}{\underline{\mu}}\cdot  t\|\nabla u\|_{L^2}^{12}  + \|\sqrt{\rho} u_t\|_{L^2}^2 \right)    \, dt  \right] \\
&\ \ \ \ \ \ \cdot \exp\left\{ \int_0^T \left[  \left( \frac{C \overline{\rho}^3 }{\underline{\mu}^3} + \frac{CM_2^2 \overline{\rho}^3}{\underline{\mu}} \right) \|\nabla u\|_{L^2}^4 +  \frac{C M^2 M_2^4 \overline{\rho}^2}{\underline{\mu}}  \|\sqrt{\rho} u_t\|_{L^2}^2      \right]\, dt  \right\}.
\ea \ee
Taking \eqref{3.2} and \eqref{3.7} into account,
\be \la{4-5-12} \ba
& \sup_{t\in [0, T] } t\int \rho |u_t|^2\, dx +\underline{\mu} \int_0^T t \|\nabla u_t\|_{L^2}^2 \\
& \leq \left[  \int_0^T  \left(  \frac{C M_2^4 \overline{\rho}^6}{\underline{\mu}}\cdot t\|\nabla u\|_{L^2}^{10} +
\frac{CM^2 M_2^8 \overline{\rho}^8}{\underline{\mu}}\cdot  t\|\nabla u\|_{L^2}^{12}  + \|\sqrt{\rho} u_t\|_{L^2}^2 \right)    \, dt  \right] \\
&\ \ \ \ \ \ \cdot \exp\left\{C\left( \frac{\overline{\rho}^4 }{\underline{\mu}^4 } + \frac{ M_2^2 \overline{\rho}^4 }{\underline{\mu}^2 } +  M^2 M_2^4 \overline{\rho}^2 \right)\right\}.
\ea \ee
According to Theorems \ref{energy-inequality-3}, \ref{corollary} and the assumption \eqref{condition2},
\be \nonumber
\ba
\int_0^T t \|\nabla u\|_{L^2}^{10} \, dt
&\leq \sup_{t\in [0, T]} t \|\nabla u\|_{L^2}^2 \cdot \sup_{t\in [0, T]} \|\nabla u\|_{L^2}^6\cdot \int_0^T \|\nabla u\|_{L^2}^2 \, dt \\
& \leq \frac{C \cdot \overline{\rho}^2}{\underline{\mu}^2} \|u_0\|_{L^2}^2 \leq \frac{C \cdot \overline{\rho}^2}{\underline{\mu}^2} \|\nabla u_0\|_{L^2}^2.\\
 \ea
 \ee
Similarly,
\be \nonumber
\int_0^T t \|\nabla u\|_{L^2}^{12} \, dt \leq \frac{C \cdot \overline{\rho}^2}{\underline{\mu}^2} \|\nabla u_0\|_{L^2}^2.
\ee
And by virture of Theorem \ref{prop1-section4},
\be \nonumber
\int_0^T \|\sqrt{\rho} u_t\|_{L^2}^2 \, dt \leq \underline{\mu} \cdot  \|\nabla u_0\|_{L^2}^2 .
\ee
Hence,
\be \la{4-5-15}  \ba
&\sup_{t\in [0, T] } t\int \rho |u_t|^2\, dx +\underline{\mu} \int_0^T t \|\nabla u_t\|_{L^2}^2 \, dt \\
&\leq C\|\nabla u_0\|_{L^2}^2\left(  \frac{ M_2^4 \overline{\rho}^8}{\underline{\mu}^3} + \frac{ M^2 M_2^8 \overline{\rho}^{10}}{\underline{\mu}^3} + \underline{\mu} \right) \cdot
\exp\left\{C\left( \frac{\overline{\rho}^4 }{\underline{\mu}^4 } + \frac{ M_2^2 \overline{\rho}^4 }{\underline{\mu}^2 } +  M^2 M_2^4 \overline{\rho}^2 \right)\right\}.
\ea\ee

On the other hand, multiplying \eqref{4-5-10} by $t$, one has
\be \la{4-5-18} \ba
 &   \frac{d}{dt} \left( \frac{ t^2}{2} \int \rho |u_t|^2 \, dx \right) +  \underline{\mu} t^2 \|\nabla u_t\|_{L^2}^2 \\
& \leq \frac{C \overline{\rho}^3 }{\underline{\mu}^3} \cdot t^2\|\sqrt{\rho} u_t\|_{L^2}^2 \cdot \|\nabla u\|_{L^2}^4
                  + \frac{C M_2^2 \overline{\rho}^3}{\underline{\mu}}  \cdot t^2\|\sqrt{\rho} u_t\|_{L^2}^2 \cdot \|\nabla u\|_{L^2}^4  \\ &\ \ \ \ + \frac{C M_2^4 \overline{\rho}^6}{\underline{\mu}} \cdot t^2\|\nabla u\|_{L^2}^{10}
+ \frac{C M^2 M_2^4 \overline{\rho}^2 }{\underline{\mu}} \cdot t^2\|\sqrt{\rho} u_t\|_{L^2}^4 \\
&\ \ \ \ + \frac{C M^2 M_2^8 \overline{\rho}^8}{\underline{\mu}}\cdot t^2 \|\nabla u\|_{L^2}^{12}
+ t \int \rho |u_t|^2 \, dx.
\ea \ee
Applying Gronwall's inequality,
\be \la{4-5-19} \ba
&\sup_{t\in [0, T] } t^2 \int \rho |u_t|^2\, dx +\underline{\mu} \int_0^T t^2 \|\nabla u_t\|_{L^2}^2 \, dt \\
& \leq \left[  \int_0^T  \left(  \frac{C M_2^4 \overline{\rho}^6}{\underline{\mu}}\cdot t^2\|\nabla u\|_{L^2}^{10} +
\frac{CM^2 M_2^8 \overline{\rho}^8}{\underline{\mu}}\cdot  t^2\|\nabla u\|_{L^2}^{12}  + t\|\sqrt{\rho} u_t\|_{L^2}^2 \right)    \, dt  \right] \\
&\ \ \ \ \ \ \cdot \exp\left\{C\left( \frac{\overline{\rho}^4 }{\underline{\mu}^4 } + \frac{ M_2^2 \overline{\rho}^4 }{\underline{\mu}^2 } +  M^2 M_2^4 \overline{\rho}^2 \right)\right\}.
\ea \ee
According to Theorems \ref{energy-inequality-3} and  \ref{corollary},
\be
\ba\nonumber
\int_0^T t^2  \|\nabla u\|_{L^2}^{10} \, dt
&\leq \sup_{t\in [0, T]} t^2 \|\nabla u\|_{L^2}^4 \cdot \sup_{t\in [0, T]} \|\nabla u\|_{L^2}^4\cdot \int_0^T \|\nabla u\|_{L^2}^2 \, dt\\
& \leq \frac{C \overline{\rho}^3}{\underline{\mu}^3} \|\nabla u_0\|_{L^2}^2.
\ea
\ee
Similarly,
\be \nonumber
\int_0^T t^2 \|\nabla u\|_{L^2}^{12} \, dt \leq \frac{C \overline{\rho}^3}{\underline{\mu}^3} \|\nabla u_0\|_{L^2}^2.
\ee
Hence,
\be \la{4-5-20}  \ba
&\sup_{t^2\in [0, T] } t^2\int \rho |u_t|^2\, dx +\underline{\mu} \int_0^T t^2 \|\nabla u_t\|_{L^2}^2  \, dt \\
&\leq C\|\nabla u_0\|_{L^2}^2\left(  \frac{ M_2^4 \overline{\rho}^9}{\underline{\mu}^4} + \frac{ M^2 M_2^8 \overline{\rho}^{11}}{\underline{\mu}^4} +  \overline{\rho} \right) \cdot
\exp\left\{C\left( \frac{\overline{\rho}^4 }{\underline{\mu}^4 } + \frac{ M_2^2 \overline{\rho}^4 }{\underline{\mu}^2 } +  M^2 M_2^4 \overline{\rho}^2 \right)\right\}.
\ea\ee

\end{proof}

%%%%%%%%%%%%%%%%%%%%%%%%%%%%%%%

\bl \la{corollary2}
Suppose $(\rho, u, P)$ is the unique local strong solution to \eqref{ns} on $[0, T]$, with the initial data $(\rho_0, u_0)$,  and satisfies
the assumptions \eqref{condition1}-\eqref{condition3}.
Then for any $r\in (3,\min\{q,6\})$
\be \la{4-6-1} \ba
& \int_0^T \|\nabla u\|_{L^\infty} \, dt \\
& \leq   C\|\nabla u_0\|_{L^2} \left[M_r \overline{\rho}^{\frac{5r-6}{4r}}\underline{\mu}^{-\frac{3(r-2)}{4r}} \left(1+\frac{\overline{\rho}}
{\underline{\mu}}\right)^{\frac12}\Theta_1^{\frac12}\exp\{C\Theta_2\}
+ M_r^{\frac{5r-6}{r}}\frac{\overline{\rho}^{\frac{6(r-1)}{r}}}
{\underline{\mu}}\right].
\ea \ee
\el

\begin{proof}
By virture of Lemma \ref{stokes}, one has for $r\in (3,\min\{q,6\})$
\be \la{4-6-2} \ba
\|\nabla u\|_{W^{1, r}} &\leq C M_r \left(  \|\rho u_t\|_{L^r} + \|\rho u\cdot \nabla u\|_{L^r}        \right)\\
& \leq C M_r \left( \|\rho u_t\|_{L^2}^{\frac{6-r}{2r}} \cdot \|\rho u_t\|_{L^6}^{\frac{3(r-2)}{2r}} +
\overline{\rho} \|u\|_{L^6} \cdot \|\nabla u\|_{L^{6r/(6-r)}}    \right)\\
& \leq C M_r  \left( \|\rho u_t\|_{L^2}^{\frac{6-r}{2r}} \cdot \|\rho u_t\|_{L^6}^{\frac{3(r-2)}{2r}} +
\overline{\rho} \|\nabla u\|_{L^2}^{\frac{6(r-1)}{5r-6}} \cdot \|\nabla u\|_{W^{1, r}}^{\frac{4r-6}{5r-6}}    \right).\\
\ea \ee
Applying Young's inequality and Sobolev embedding inequality,
\be \la{4-6-3}
\|\nabla u\|_{W^{1, r}} \leq C M_r \overline{\rho}^{\frac{5r-6}{4r}} \cdot \|\sqrt{\rho} u_t\|_{L^2}^{\frac{6-r}{2r}}
\cdot \|\nabla u_t\|_{L^2}^{\frac{3(r-2)}{2r}} + CM_r^{\frac{5r-6}{r}} \overline{\rho}^{\frac{5r-6}{r}}\cdot \|\nabla u\|_{L^2}^{\frac{6(r-1)}{r}}.
\ee
Hence,
\be \la{4-6-4} \ba
& \int_0^T \|\nabla u\|_{L^\infty} \, dt \\ & \leq C \int_0^T \|\nabla u\|_{W^{1, r}}\, dt \\
& \leq   C\int_0^T \left( M_r \overline{\rho}^{\frac{5r-6}{4r}}  \|\sqrt{\rho} u_t\|_{L^2}^{\frac{6-r}{2r}}
 \|\nabla u_t\|_{L^2}^{\frac{3(r-2)}{2r}} + M_r^{\frac{5r-6}{r}} \overline{\rho}^{\frac{5r-6}{r}} \|\nabla u\|_{L^2}^{\frac{6(r-1)}{r}} \right) \, dt.
\ea \ee

If $T\leq 1$, according to Theorem \ref{prop2-section4},
\be \nonumber \ba
& \int_0^T \|\sqrt{\rho} u_t\|_{L^2}^{\frac{6-r}{2r}} \cdot  \|\nabla u_t\|_{L^2}^{\frac{3(r-2)}{2r}} \, dt \\
& \leq \int_0^T \left(  t^{\frac12} \|\sqrt{\rho} u_t\|_{L^2}\right)^{\frac{6-r}{2r}} \cdot \left(t^{\frac12} \|\nabla u_t\|_{L^2}  \right)^{\frac{3(r-2)}{2r} } \cdot t^{-\frac12} \, dt \\
& \leq   C \left( \sup_{t\in [0, T]} t^{\frac12} \|\sqrt{\rho} u_t\|_{L^2}\right)^{\frac{6-r}{2r} }
\cdot \left( \int_0^T t\|\nabla u_t \|_{L^2}^2\, dt \right)^{\frac{3(r-2)}{4r}}\left(\int_0^Tt^{-\frac{2r}{r+6}}dt
\right)^{\frac{r+6}{4r}} \\
& \leq  C\underline{\mu}^{-\frac{3(r-2)}{4r}} \cdot \|\nabla u_0\|_{L^2}^2 \Theta_1^{\frac12}\exp\{C\Theta_2\}.
\ea \ee

If $T>1$, applying Theorem \ref{prop2-section4} again,
\be \nonumber \ba
& \int_0^T \|\sqrt{\rho} u_t\|_{L^2}^{\frac{6-r}{2r}} \cdot  \|\nabla u_t\|_{L^2}^{\frac{3(r-2)}{2r}} \, dt \\
& \leq \int_0^T \left(  t^{\frac12} \|\sqrt{\rho} u_t\|_{L^2}\right)^{\frac{6-r}{2r}} \cdot \left(t^{\frac12} \|\nabla u_t\|_{L^2}  \right)^{\frac{3(r-2)}{2r} } \cdot t^{-\frac12} \, dt \\
&\ \ \ \ + \int_1^T \left(t \|\sqrt{\rho  } u_t\|_{L^2} \right)^{\frac{6-r}{2r}} \cdot  \left(t  \|\nabla u_t\|_{L^2}  \right)^{\frac{3(r-2)}{2r} } \cdot t^{-1} \, dt\\
& \leq C \left( \sup_{t\in [0, T]} t^{\frac12} \|\sqrt{\rho} u_t\|_{L^2}\right)^{\frac{6-r}{2r} }
\cdot \left( \int_0^T t\|\nabla u_t \|_{L^2}^2\, dt \right)^{\frac{3(r-2)}{4r}}\left(\int_0^Tt^{-\frac{2r}{r+6}}dt
\right)^{\frac{r+6}{4r}} \\
&\ \ \ \ + C \left( \sup_{t \in [0, T]}   t\|\sqrt{\rho} u_t\|_{L^2}\right)^{\frac{6-r}{2r} }
\cdot \left( \int_0^T t^2 \|\nabla u_t \|_{L^2}^2\, dt \right)^{\frac{3(r-2)}{4r}}\left(\int_1^Tt^{-\frac{4r}{r+6}}dt\right)^{\frac{6+r}{4r}}\\
& \leq C\underline{\mu}^{-\frac{3(r-2)}{4r}} \|\nabla u_0\|_{L^2} \Theta_1^{\frac12}\exp\{C\Theta_2\} + C \underline{\mu}^{-\frac{3(r-2)}{4r}} \left(\frac{\overline{\rho}}{\underline{\mu}}\right)^{\frac12}\|\nabla u_0\|_{L^2} \Theta_1^{\frac12}\exp\{C\Theta_2\}\\
& \leq C \underline{\mu}^{-\frac{3(r-2)}{4r}} \left(1+\frac{\overline{\rho}}{\underline{\mu}}\right)^{\frac12}\|\nabla u_0\|_{L^2} \Theta_1^{\frac12}\exp\{C\Theta_2\}.
\ea \ee
On the other hand,
\be
\int_0^T \|\nabla u\|_{L^2}^{\frac{6(r-1)}{r}} \, dt  \leq \sup_{t\in [0, T]} \|\nabla u\|_{L^2}^{\frac{4r-6}{r}} \cdot \int_0^T \|\nabla u\|_{L^2}^2 \, dt \leq \frac{C \cdot \overline{\rho}}{\underline{\mu}} \|\nabla u_0\|_{L^2}.
\ee

Therefore,
%\be \nonumber \ba
%& \int_0^T \|\nabla u\|_{L^\infty} \, dt  \\
%& \leq \left[  CM_4 \overline{\rho}^{\frac78}   \left(  \frac{C M_2^2 \overline{\rho}^4}{\underline{\mu}^{\frac32}} + \frac{C M M_2^4 \overline{\rho}^{5}}{\underline{\mu}^{\frac32}} + C\underline{\mu}^{\frac12} \right) \cdot
%\exp\left\{ \frac{C\overline{\rho}^4 }{\underline{\mu}^4 } + \frac{C M_2^2 \overline{\rho}^4 }{\underline{\mu}^2 } + C M^2 M_2^4 \overline{\rho}^2 \right\}  \right. \\
%&\ \ \ \ +CM_4 \overline{\rho}^{\frac78} \left(  \frac{C M_2^4 \overline{\rho}^{\frac92} }{\underline{\mu}^2} + \frac{C M  M_2^4 \overline{\rho}^{\frac{11}{2} }}{\underline{\mu}^2} + C \overline{\rho}^{\frac12} \right) \cdot
%\exp\left\{ \frac{C\overline{\rho}^4 }{\underline{\mu}^4 } + \frac{C M_2^2 \overline{\rho}^4 }{\underline{\mu}^2 } + C M^2 M_2^4 \overline{\rho}^2 \right\}\\
%&\ \ \ \ \ \left.   + \frac{CM_4^{\frac72} \overline{\rho}^{\frac92}}{\underline{\mu}}     \right]\cdot \|\nabla u_0\|_{L^2}\\
%& \triangleq C_2(M,\overline{\rho},\underline{\mu},\overline{\mu}) \|\nabla u_0\|_{L^2}.
%\ea \ee
\be \la{4-6-1} \ba
& \int_0^T \|\nabla u\|_{L^\infty} \, dt \\
& \leq   C\|\nabla u_0\|_{L^2} \left[M_r \overline{\rho}^{\frac{5r-6}{4r}}\underline{\mu}^{-\frac{3(r-2)}{4r}} \left(1+\frac{\overline{\rho}}
{\underline{\mu}}\right)^{\frac12}\Theta_1^{\frac12}\exp\{C\Theta_2\}
+ M_r^{\frac{5r-6}{r}}\frac{\overline{\rho}^{\frac{6(r-1)}{r}}}
{\underline{\mu}}\right].\\
& \triangleq C_2(M,\overline{\rho},\underline{\mu},\overline{\mu}) \|\nabla u_0\|_{L^2}.
\ea \ee

\end{proof}

%%%%%%%%%%%%%%%%%%%%%%%%%%%%%%%%
\bt \la{prop3-section4} Suppose $(\rho, u, P)$ is the unique local strong solution to \eqref{ns} on $[0, T]$ and
the assumptions \eqref{condition1}-\eqref{condition3}.
There exists a positive number $\epsilon_0$, depending on $\Omega$, $q$, $M$, $\overline{\rho}$, $\overline{\mu}$, $\underline{\mu}$ such that if
\be \nonumber \|\nabla u_0\|_{L^2} \leq \epsilon_0~,
\ee
then
\be \la{4-8-1}
\sup_{t\in [0, T]} \|\nabla \mu(\rho)\|_{L^q} \leq 2M~,
\ee
and
\be \nonumber
\sup_{t\in [0, T]} \|\nabla \rho\|_{L^q} \leq 2\|\nabla\rho_0\|_{L^q}.
\ee

\et

\begin{proof}
Consider the $x_i$-derivative of the renormalized mass equation for $\mu(\rho)$,
\be \nonumber
\left ( \partial_i \mu(\rho)\right)_t + (\partial_i u \cdot \nabla)  \mu(\rho) + u \cdot \nabla \partial_i \mu(\rho) =0~.
\ee
It implies that for every $t\in [0, T]$, in view of Lemma \ref{corollary2}, one has
\be \la{4-7-1}  \ba
\|\nabla \mu(\rho)(t)  \|_{L^q}  &  \leq \|\nabla \mu(\rho_0)\|_{L^q}  \cdot \exp\left\{\int_0^t \|\nabla u\|_{L^\infty} \, ds\right\}\\
& \leq \|\nabla \mu(\rho_0)\|_{L^q}  \cdot \exp\left\{ C_2(M,\overline{\rho},\underline{\mu},
\overline{\mu})  \cdot \|\nabla u_0\|_{L^2}\right\}.
\ea
\ee
Choose some small positive constant $\epsilon_0$, which satisfies
\be\nonumber
\epsilon_0 \leq \frac12, \ \ \ \ \epsilon_0^2 \cdot C_1 \underline{\mu}^{-2} (M_2^2 + M^4 M_2^6) \overline{\rho}^4
\leq \ln 2,
\ee
and
\be \nonumber
\epsilon_0 \cdot C_2(M, \overline{\rho}, \overline{\mu}, \underline{\mu}) \leq \ln 2.						
\ee
If $\|\nabla u_0\|_{L^2} \leq \epsilon_0$, \eqref{4-8-1} holds.

Similarly,
\be \la{4-8-3} \ba
\|\nabla \rho(t)  \|_{L^q}  & \leq \|\nabla \rho_0\|_{L^q}  \cdot \exp\left\{\int_0^t \|\nabla u\|_{L^\infty} \, ds\right\}\\
& \leq 2\|\nabla \rho_0\|_{L^q},
\ea
\ee
Therefore,  Theorem \ref{prop3-section4} is proved.
\end{proof}

%%%%%%%%%%%%%%Subsection 3.2%%%%%%%%%%%%%%%%%%%%%%%%%%%%%
\subsection{Proof of Theorem \ref{main-result}}
With the a priori estimates in Subsection 3.1 in hand, we are prepared for the proof of Theorem \ref{main-result}.
\begin{proof}
According to Theorem \ref{main-result-CK}, there exists a $T_*>0$ such that the density-dependent Navier-Stokes system \eqref{ns} has a unique local strong solution $(\rho, u, P)$ on $[0, T_*]$. We plan to extend the local solution to a global one.

Since $\|\nabla \mu(\rho_0)\|_{L^q}=M <4M$, and due to the continuity of $\nabla \mu(\rho)$ in $L^q$ and $\nabla u_0$ in $L^2$,  there exists a $T_1 \in (0, T_*)$ such that
$\sup_{0\leq t \leq T_1}  \|\nabla \mu(\rho)(t)\|_{L^q} \leq 4M, $
and at the same time
$\sup_{0\leq t \leq T_1} \|\nabla u(t)\|_{L^2} \leq 2\|\nabla u_0\|_{L^2}.$

Set
\be \nonumber T^* = \sup \left\{T | \ (\rho, u, P)\  \mbox{is a strong solution to} \ \eqref{ns}\ \mbox{on} \ [0, T]   \right\}~.
\ee
\be \nonumber
  T_1^* = \sup \left\{   T \big|   \ba  & \ (\rho, u, P)\   \mbox{is a strong solution to} \ \eqref{ns}\ \mbox{on} \ [0, T], \\ & \sup_{0 \leq t \leq T}  \|\nabla \mu(\rho) \|_{L^q} \leq 4M ,  \ \mbox{and}\ \sup_{0 \leq t \leq T}   \|\nabla u\|_{L^2} \leq 2 \|\nabla u_0\|_{L^2}.
\ea  \right\}.
 \ee
Then $T_1^* \geq T_1 >0$. Recalling Theorems \ref{prop1-section4} and \ref{prop3-section4}, it's easy to verify
\be
T^* = T_1^*,
\ee
provided that $\|\nabla u_0\|_{L^2} \leq \epsilon_0$ as assumed.

We claim that $T^* = \infty$. Otherwise,  assuming that $T^* < \infty$. By virture of Theorems \ref{prop1-section4} and  \ref{prop3-section4}, for every $t\in [0, T^*)$, it holds that
\be \label{E1}
 \|\nabla \rho (t) \|_{L^q} \leq 2\|\nabla \rho_0\|_{L^q}, \ \ \mbox{and}\ \|\nabla u(t)\|_{L^2} \leq \sqrt{2} \|\nabla u_0\|_{L^2}.
\ee
which contradicts to the blowup criterion (\ref{criterion-CK}). Hence we complete  the proof for Theorem \ref{main-result}.

\end{proof}

{\bf Acknowledgement:}\ The research of Xiangdi Huang is supported in part by NSFC Grant No. 11101392. The research of Yun Wang is supported in part by NSFC Grant No. 11301365.

\end{document}